\documentclass[a4paper,12pt]{amsart}
\usepackage[latin1]{inputenc}
\usepackage[T1]{fontenc}
\usepackage{amsfonts,amssymb}
\usepackage{pstricks}
\usepackage{amsthm}
\usepackage[all]{xy}
\usepackage[pdftex]{graphicx}
\usepackage{graphicx}   
\usepackage{xcolor}
\newtheorem{theo}{Theorem}[section]

\newtheorem{rem}[theo]{Remark}

\title[Even-Odd partition identities of G\"{o}llnitz--Gordon type]{Even-Odd partition identities of G\"{o}llnitz--Gordon type}

\author{Pooneh Afsharijoo}

\begin {document}
\maketitle
\begin{abstract} In this paper, we establish a new companion to the G\"ollnitz--Gordon identities. We give two independent proofs. The first is based on a new system of recurrence formulas, while the second uses generating series and the analytic form of the G\"ollnitz--Gordon identities.
\end{abstract}
\section{Introduction}
A partition of a non-negative integer $n$ is a non-increasing  sequence $\lambda=(\lambda_1\geq \cdots \geq \lambda_m)$ of positive integers $\lambda_i$ (called parts of $\lambda$) whose sum is  $n$. 
\\
A partition identity is an equality between the number of partitions of an integer $n$ satisfying a condition $A$ and the number of partitions of $n$ satisfying another condition $B$ which is true for any  $n\geq0.$ Such identities play an important role in many areas such as number theory, statistical physics and Lie theory. 
\bigskip

 Let $n$ be a positive integer. For $i \in \{1,2\}$ let  $\mathcal{G}_{i}(n)$ denote the set of partitions of $n$ whose parts are at least $2i-1$, with no equal parts, no consecutive parts, and no consecutive even parts. Let $\mathcal{G}'_{i}(n)$ denote the set of partitions of $n$ whose parts are congruent to $4, (2i-1)$ or $(9-2i)$ modulo $8$.
We have the following classical theorem.
\begin{theo}(G\"{o}llnitz--Gordon identities) For the integers $n\geq 0$ and $i\in\{1,2\}$ let $G_{i}(n)$  and $G'_{i}(n)$ denote the cardinalities of $\mathcal{G}_{i}(n)$ and $\mathcal{G}'_{i}(n)$ respectively. Then we have:
 $$G_{i}(n)=G'_{i}(n).$$
\end{theo}
For example,
$$\mathcal{G}_{1}(10)=\{10, 9+1, 8+2, 7+3, 6+3+1\},$$
$$\mathcal{G}'_{1}(10)=\{9+1, 7+1+1+1, 4+4+1+1, 4+1+1+1+1+1+1, 1+1+1+1+1+1+1+1+1+1\},$$
$$\mathcal{G}_{2}(10)=\{10, 7+3\},$$
$$\mathcal{G}'_{2}(10)=\{5+5, 4+3+3\}.$$
Thus $G_1(10)=G'_1(10)=5$ and $G_2(10)=G'_2(10)=2.$
\bigskip

In this paper, we establish a new companion to the G\"ollnitz--Gordon identities (Theorem~\ref{k=1}). This new companion treats the even and odd parts of a partition differently. To state the theorem, we first introduce some notation. Let $\lambda$ be a partition of an integer $n$. 
\bigskip

We denote by $\ell(\lambda)$ the number of parts of $\lambda,$ and by $|\lambda|$ the sum of its parts. We write $\lambda=(\lambda^{(e)},\lambda^{(o)}),$ where $\lambda^{(e)}$ and $\lambda^{(o)}$ denote the subpartitions consisting of the even and odd parts of $\lambda$, respectively.
We now state our main result.
\begin{theo}\label{k=1}

 (A new companion to the G\'ollnitz--Gordon identities)

  Given integers $n\ge0$ and $i\in\{1,2\}$, let $\mathcal{H}_{i}(n)$
 be the set of partitions $\lambda$ of $n$ whose odd parts are distinct and at least  $2i-1,$ and whose even parts are at least $2(\ell(\lambda)+i-1).$ Let $H_i(n)$ denote the cardinality of $\mathcal{H}_{i}(n)$. Then
$H_i(n)=G_i(n)=G'_{i}(n)$.



\end{theo}
For example,
$$\mathcal{H}_1(10)=\{10, 9+1, 7+3, 6+4, 6+3+1\},$$
and 
$$\mathcal{H}_2(10)=\{10, 7+3\}.$$
Hence, $G_1(10)=G'_1(10)=H_1(10)=5$ and $G_2(10)=G'_2(10)=H_2(10)=2.$
\bigskip

 In \cite{Af2}, we developed a method based on recursion formulas for proving partition identities, which yielded new companions to the Rogers--Ramanujan identities. In Section~\ref{P1}, we adapt this method to prove Theorem~\ref{k=1}.  Let \[
\mathcal{H}_i(r,s,n)=
\{
\lambda\in\mathcal{H}_i(n):
\ell(\lambda^{(e)})=r,\;
\ell(\lambda^{(o)})=s
\}.
\] 
Similarly, define \[
\mathcal{G}_{i}(r,s,n)=
\{
\lambda\in\mathcal{G}_i(n):
\ell(\lambda^{(e)})=r,\;
\ell(\lambda^{(o)})=s
\}.
\] 
Let $H_i(r,s,n)$ (respectively $G_i(r,s,n)$) denote the cardinality of $\mathcal{H}_{i}(r,s,n)$ (respectively  $\mathcal{G}_{i}(r,s,n)$). We  \textit{introduce} a new system of recursive formulas between the $H_i(r,s,n)$s. Next we prove that the $G_i(r,s,n)$s also satisfy the same system. This completes the proof since the $H_i(r,s,n)$s are uniquely determined by this system. 
\bigskip

In Section \ref{P2}, we give a second proof of Theorem~\ref{k=1} using the analytic form of the G\"ollnitz--Gordon identities. We show that the generating series of $H_i(n)$ coincides with that of $G_i(n).$



\section{First proof of Theorem \ref{k=1}: recursion formulas }\label{P1}
In this section, we give a first proof of Theorem~\ref{k=1}.

 \begin{proof}(Proof of the Theorem \ref{k=1})
 It suffices to show that both  $H_{i}(r,s,n)$ and $G_{i}(r,s,n)$  satisfy the following system of recursion formulas.

 \begin{equation}\label{1}
 \begin{split}
 &H_{i}(r,s,n)=\begin{cases}
1 &\text{if } r=s=n=0  \\
0 &\text{if }\min\{r,s,n\}\leq0,\ (r,s,n)\neq(0,0,0);
\end{cases}\\
\end{split}
\end{equation} 
\begin{equation}\label{2}
H_{1}(r,s,n)-H_{2}(r,s,n)=H_{1}(r,s-1,n-2(r+s)+1)+H_{1}(r-1,s,n-4(r+s)+2); 
\end{equation}

\begin{equation}\label{3}
H_{2}(r,s,n)=H_{1}(r,s,n-2(r+s)).
\end{equation}

The initial conditions are immediate. Indeed, 0 has a unique partition (the empty partition), whereas a negative integer has no partitions. Likewise, a positive integer has no partition of non-positive length.
\bigskip

In order to prove that $H_{i}(r,s,n)$ satisfy equation (\ref{2}), note that $\mathcal{H}_{1}(r,s,n)\setminus \mathcal{H}_{2}(r,s,n)$ contains the partitions of $n$ of the form  $$\lambda:(\underbrace{\lambda_{i_1} \geq \cdots \geq \lambda_{i_{r}}}_{\text {The even parts of $\lambda$}},\underbrace{\lambda_{j_1}> \cdots >\lambda_{j_{s}}}_{\text {The odd parts of $\lambda$}}),$$

with either $\lambda_{j_s}=1$ or $\lambda_{i_r}=2(r+s).$  
\\
\begin{itemize}
\item \textbf{Case $1$}: Suppose that $\lambda_{j_s}=1.$ In this case we define a bijection $\phi$ by removing a bijection $\phi$ to a partition $\mu$ by removing $\lambda_{j_s}=1$ and then subtracting $2$ from each other part of $\lambda.$ We obtain:
 $$\mu:(\underbrace{\lambda_{i_1}-2 \geq \cdots \geq \lambda_{i_{r}}-2}_{\text {The even parts of $\mu$}},\underbrace{\lambda_{j_1}-2> \cdots > \lambda_{j_{s-1}}-2}_{\text {The odd parts of $\mu$}}),$$
 
 which is a partition of $n-2(r+s)+1$ with exactly $r$ even parts, each  at least equal to $2(r+s)-2$, and $s-1$ distinct odd parts at least equal to $1.$ Therefore $\mu\in \mathcal{H}_{1}(r,s-1,n-2(r+s)+1).$
 \\
 \item \textbf{Case $2$}: Suppose that  $\lambda_{j_s} \neq1.$ Then necessarily $\lambda_{i_r}=2(r+s).$ This time, we transform the partition  $\lambda$ to a partition $\mu$ by  removing $\lambda_{i_r}$ and then subtracting $2$ from each other part. We obtain:
 $$\mu:(\underbrace{\lambda_{i_1}-2\geq \cdots \geq \lambda_{i_{r-1}}-2}_{\text {The even parts of $\mu$}},\underbrace{\lambda_{j_1}-2>\cdots > \lambda_{j_{s}}-2}_{\text {The odd parts of $\mu$}})\in \mathcal{H}_{1}(r-1,s,n-4(r+s)+2).$$
 
\end{itemize}

By construction, both transformations are bijections. Therefore, $H_i(r,s,n)$ satisfy equation (\ref{2}).

 In order to prove this equation for $G_{i}(r,s,n)$, note that, $\mathcal{G}_{1}(r,s,n) \setminus \mathcal{G}_{2}(r,s,n)$ is the set of partitions of $n$ of the form:
  $$\lambda:(\underbrace{\lambda_{i_1} > \cdots > \lambda_{i_{r}}}_{\text {The even parts of $\lambda$}}, \underbrace{\lambda_{j_1} > \cdots > \lambda_{j_{s}}}_{\text {The odd parts of $\lambda$}}),$$

Without consecutive parts, without consecutive even parts and either $\lambda_{j_s}=1$ or $\lambda_{i_r}=2.$
\\
\begin{itemize}
\item \textbf{Case $1$}: Suppose that $\lambda_{j_r}=1.$ Since the parts of $\lambda$ are not equal or consecutive, we have $\lambda_{j_{s-1}}\geq 3$ and $\lambda_{i_r}\geq 4.$  We transform $\lambda$ to a partition $\mu$ by the bijection $\phi$ defined above. It is straightforward to verify that $\mu$ belongs to $ \mathcal{G}_{1}(r,s-1,n-2(r+s)+1).$ 
\\
\item \textbf{Case $2$}: Suppose that $\lambda_{j_r} \neq 1.$ Then $\lambda_{i_r}=2.$ Note that since $\lambda$ does not have equal or consecutive even parts nor equal or consecutive  parts, so $\lambda_{i_{r-1}}\geq 6$ and $\lambda_{j_s}\geq 5.$ Now we remove $\lambda_{i_r}$ and then we subtract $4$ from each other part. We obtain:
$$\mu:(\underbrace{\lambda_{i_1}-4>\cdots> \lambda_{i_{r-1}}-4}_{\text {The even parts of $\mu$}},\underbrace{\lambda_{j_1}-4>\cdots > \lambda_{j_{s}}-4}_{\text {The odd parts of $\mu$}}).$$
It is immediate that $\mu$ belongs to $\mathcal{G}_{1}(r-1,s,n-4(r+s)+2).$

\end{itemize}
The bijectivity of these transformations proves that $G_{i}(r,s,n)$ satisfy equation (\ref{2}). We next prove (\ref{3}) for both $H_{i}(r,s,n)$ and $G_{i}(r,s,n).$ Let $$\lambda:(\underbrace{\lambda_{i_1} \geq \cdots \geq \lambda_{i_{r}}}_{\text {The even parts of $\lambda$}},\underbrace{\lambda_{j_1}>\cdots> \lambda_{j_{s}}}_{\text {The odd parts of $\lambda$}})\in \mathcal{H}_{2}(r,s,n).$$ Apply the bijection $\psi$ on $\lambda$, which subtracts $2$ from every part we obtain a partition $\mu\in \mathcal{H}_{1}(r,s,n-2(r+s)),$ and  (\ref{3}) follows for $H_{i}(r,s,n).$
\\
In order to prove this last equation for $G_{i}(r,s,n)$ we send a partition $$\lambda:(\underbrace{\lambda_{i_1}> \cdots>\lambda_{i_{r}}}_{\text {The even parts of $\lambda$}},\underbrace{\lambda_{j_1}> \cdots > \lambda_{j_{s}}}_{\text {The odd parts of $\lambda$}})\in \mathcal{G}_{2}(r,s,n)$$ to a partition $\mu$ by the bijection $\psi$. 
Since $\lambda_{i_r}\geq4$ and $\lambda_{j_s}\geq 3$ thus the smallest even part of $\mu$ (respectively the smallest odd part of $\mu$) is at least equal to $2$ (respectively is at least equal to $1$). Therefore, $\mu\in \mathcal{G}_{1}(r,s,n-2(r+s)).$  Thus, (\ref{3}) also holds for$G_{i}(r,s,n).$
\\
We have now shown that both $H_{i}(r,s,n)$ and $G_{i}(r,s,n)$  satisfy (\ref{1})--(\ref{3}). Finally, the recurrence system together with the initial conditions uniquely determines the numbers $H_i(r,s,n)$ by a triple induction on $r,s$ and $n$. Hence $H_{i}(r,s,n)=G_{i}(r,s,n)$ for all integers $r,s,n\geq 0$ and  $i\in\{1,2\}.$ Summing over all r,s completes the proof of Theorem~\ref{k=1}. 

\end{proof}

\begin{rem} The proof above establishes a refinement of Theorem~\ref{k=1}. More precisely, for all integers $r,s,n\ge0$ and $i\in\{1,2\}$,
\[
H_i(r,s,n)=G_i(r,s,n).
\]
In other words, the identity remains valid after fixing the numbers of even and odd parts.
\end{rem}

\section{Second proof of Theorem \ref{k=1}: Generating series}\label{P2}

In this section, we give a second proof of Theorem~\ref{k=1} using generating series. We first recall the definition of the $q$-binomial coefficient: \[
\begin{bmatrix}
n\\
k
\end{bmatrix}_q
=
\frac{(q;q)_n}{(q;q)_r\,(q;q)_{n-r}},
\qquad 0\le k\le n.
\]
where $(a,q)_n$ denotes the $q$-Pochhammer symbol,  \[
(a;q)_n=
\begin{cases}
1, & n=0,\\[0.3em]
\displaystyle\prod_{k=0}^{n-1}(1-aq^k), & n\ge1,
\end{cases}
\qquad
(a;q)_\infty=\prod_{k=0}^{\infty}(1-aq^k).
\]
We also recall the analytic form of the G\"ollnitz--Gordon identities:
 \begin{theo}\label{A}(G\"ollnitz--Gordon identities: analytic
version)

For $i \in \{1,2\}$ we have:
\[
\sum_{n=0}^{\infty}\frac{q^{n^2+(2i-2)n}(-q,q^2)_n}{(q^2,q^2)_n}=\frac{1}{(q,q^8)_{\infty}(q^{2i-1},q^8)_{\infty}(q^{9-2i},q^8)_{\infty}}.\] 

\end{theo}
To prove Theorem~\ref{k=1}, it is enough to show that the generating series of $H_i(n)$ coincides with the left-hand side of Theorem~\ref{A}.
\begin{theo}\label{AH} For $i \in \{1,2\}$ we have:
\[
\sum_{n=0}^{\infty} H_i(n)q^n=\sum_{n=0}^{\infty}\frac{q^{n^2+(2i-2)n}(-q,q^2)_n}{(q^2,q^2)_n}.
\]

\end{theo}
\begin{proof}
Let $\lambda=(\lambda^{(e)},\lambda^{(o)})\in \mathcal{H}_i(n)$ with $r$ even parts and $s$ odd parts for some $r,s\geq 0$. The idea is to compute separately the generating seires of the even and odd subpartitions and then multiply them together. We begin with $\lambda^{(e)},$ and then that of $\lambda^{(o)}$.
\\
Suppose that $\lambda^{(e)}=(2k_1\geq \cdots \geq 2k_r).$ Since $\lambda \in \mathcal{H}_i(n),$ we have $k_r\geq r+s+i-1$. Thus, the partition
\[
\mu=(k_1\ge\cdots\ge k_r)
\]
is an ordinary partition with exactly $r$ parts, each at least
$r+s+i-1$.

Its Young diagram therefore consists of a rectangle of size
$r\times(r+s+i-1)$ together with an ordinary partition having at most
$r$ parts. Therefore, its generating series is equal to \[\frac{q^{r(r+s+i-1)}}{(q,q)_r},\] and thus the one of $\lambda^{(e)}$ is equal to \[\frac{q^{2r(r+s+i-1)}}{(q^2,q^2)_r}.\] 
\\
Let us now compute the generating series of $\lambda^{(o)}.$ Suppose that $\lambda^{(o)}=(2k'_1-1> \cdots > 2k'_s-1)$ with $k'_s\geq i.$ Thus, the partition $\mu'=(k'_1> \cdots > k'_s)$ is an ordinary partition with $s$ distinct parts, each at least equal to $i.$ In order to compute the generating series of $\mu'$. Observe that $\mu'$ can be decomposed into the staircase partition
\[
(s-1>s-2>\cdots>1)
\]
and the partition
\[
\gamma=(k'_1-(s-1)\ge k'_2-(s-2)\ge\cdots\ge k'_s).
\] The Young diagram of $\gamma$ contains a rectangle of size
$s\times i$, together with an ordinary partition having at most $s$
parts. Thus the generating series of $\mu'$ is equal to:
\[q^{s(s-1)/2}\frac{q^{si}}{(q)_s}.\] Since $|\lambda^{(o)}|=2|\mu'|-s$ it follows that the generating series of $\lambda^{(o)}$ is equal to\[q^{-s}q^{s(s-1)}\frac{q^{2si}}{(q^2)_s}=\frac{q^{s^2+2s(i-1)}}{(q^2,q^2)_s}.\] 
\\
Multiplying the generating series of $\lambda^{(e)}$ and $\lambda^{(o)}$ together and then substituting $r+s$ by $n$ we obtain:

\begin{equation}\label{series}
\begin{aligned}
\sum_{n=0}^{\infty} H_i(n)q^n=\sum_{r,s\geq 0} \Big(\frac{q^{2r(r+s+i-1)}}{(q^2,q^2)_r}\Big) \Big(\frac{q^{s^2+2s(i-1)}}{(q^2,q^2)_s}\Big)\\ =\sum_{n=0}^{\infty}\sum_{r=0}^{n}\frac{q^{n^2+r^2+2n(i-1)}}{(q^2,q^2)_r (q^2,q^2)_{n-r}} \\ =\sum_{n=0}^{\infty}\frac{q^{n^2+2n(i-1)}}{(q^2,q^2)_{n}}\sum_{r=0}^{n} \begin{bmatrix}
n\\
r
\end{bmatrix}_{q^2} q^{r^2}.
    \\ 
\end{aligned}
\end{equation}	
Let us now recall the finite $q$-binomial theorem:
\[
(-z,q)_n=\prod_{r=0}^{n-1}(1+zq^r)
=
\sum_{r=0}^{n}
q^{\binom{r}{2}}
\begin{bmatrix}
n\\
r
\end{bmatrix}_q
z^r.
\]
Applying this theorem to the inner sum in
\eqref{series}, after replacing $q$ by $q^2$ and setting $z=q$, gives
\[
\sum_{r=0}^n
\begin{bmatrix}
n\\ r
\end{bmatrix}_{q^2}
q^{r^2}
=
(-q,q^2)_n,
\]
which completes the proof.
\end{proof}
\begin{rem} Theorem~\ref{AH}, together with the analytic form of the G\"ollnitz--Gordon identities (Theorem~\ref{A}), immediately yields Theorem~\ref{k=1}.
\end{rem}
\section*{ Acknowledgment}
 I would like to thank Matthew Russell, who shared his conjecture on the new companion to the G\"ollnitz--Gordon identities presented in Theorem \ref{k=1}.
 

\bibliographystyle{alpha}
\nocite{A, Af, Af2, AB, AM, ADJM, G, GP, LZ, BMS, BMS1, M}
\bibliography{article}

\end{document}